\documentclass[12pt]{article}

\usepackage{amsmath}
\usepackage{amsfonts}     
\usepackage{amstext}
\usepackage{graphicx}
\usepackage{latexsym}
\usepackage[english]{babel}

\newtheorem{theorem}{Theorem}
\newtheorem{lemma}[theorem]{Lemma}
\newtheorem{proposition}[theorem]{Proposition}

\newtheorem{remark}[theorem]{Remark}

\begin{document}

\title{ Dimension via Waiting time and Recurrence}
\author{Stefano Galatolo}
\maketitle

\begin{abstract} 
Quantitative recurrence indicators are defined by measuring the first entrance time of 
the orbit of a point $x$ in a decreasing sequence of neighborhoods of another point $y$. It is proved that these recurrence indicators are a.e. greater or equal to the local dimension at $y$, then 
these recurrence indicators can be used to have a numerical upper bound on the local dimension of an invariant 
measure. 
\end{abstract}

\def\stackunder#1#2{\mathrel{\mathop{#2}\limits_{#1}}}

\section{Introduction} 
 
The first well known results  about recurrence  in a dynamical system $(X,T)$ state that under suitable assumptions
 a typical trajectory of the system comes back infinitely many times in any neighborhood of  its starting point. 
 These quantitative results does not give a quantitative estimation about the speed of this coming back to the starting point.

A more precise analysis of recurrence was done by defining quantitative recurrence  indicators. In the literature such indicators have been defined in several ways by measuring the first return time of an orbit in a decreasing sequence of neighborhoods of the starting point. 
These sequences  of  neighborhoods have been defined by the metric of the space $X$, considering a decreasing sequence of balls  (\cite{Bo},\cite{BS}) or with respect to the symbolic 
dynamics induced by a partition, considering a decreasing sequence of cylinders $c_k$ on the associated symbolic space  (\cite{OW}). 
Other definitions consider the forward images of the whole cylinder $c_k$ and consider as a first return for the cylinder the minimum $n$ such that $T^{n}(c_k)\cap c_k \neq \emptyset$ 
(\cite{HSV},\cite{ACS},\cite{STV},\cite{BGI},\cite{CD}). 
In the above cited papers many relations have been then  proved between these indicators and other important features of dynamics (for example dimension, entropy, orbit complexity, Lyapunov exponents, mixing properties). 
 
Barreira and Saussol in \cite{BS} prove some  relations between quantitative recurrence and the local dimension of an invariant measure. In an applicative,  framework, this relation can be used to estimate efficiently this local dimension. 
If  some technical assumptions are satisfied their recurrence indicator is indeed a.e. equal to the dimension of the invariant measure. 
Their assumptions are satisfied for example in hyperbolic systems, with an equilibrum measure  supported on a locally maximal hyperbolic set.  
 
In computer simulations or experimental situations, the quantitative recurrence indicators can be  easily estimated by looking  
to the behavior of a "typical, random" orbit and to its first entrance time in a sequence of balls centered in the starting point. 
 By the above results this can give a numerical estimation of the pointwise dimension of the underlying invariant 
measure, which is not easy to be known in general.  In general systems (where the additional assumptions are not satisfied) however the Barreira and Saussol recurrence indicator gives 
 only a {\em lower } bound on the dimension.

A natural generalization of the quantitative recurrence indicators can be defined by measuring how fast the orbit of a point $x$ approaches near to another point $y$. 
In  the literature about finite alphabet stochastic processes and symbolic systems indicators of this type were called {\em waiting times}. Relations between waiting time and entropy similar to the Oerstein Weiss theorem were proved (\cite{Sh}). Such  relations holds for Markov chains and in a weaker version in  weak Bernoulli processes. 
 In the case of finite type shifts with a Gibbs measure other results were given by \cite{Ch}. Another recent work (\cite{KS}) calculates  the waiting time
for irrational translations on the circle.

In this work we consider general systems acting on a metric space and 
we consider a recurrence-waiting time indicator that generalizes the former
ones. Then we prove some result relating this indicator to the local
dimension (of a measure defined on the space were the dynamics acts).  
 
In \cite{BGI} a quantitative recurrence argument of this kind was used  to estimate the initial condition sensitivity  
and the orbit complexity of interval exchange transformations \footnote{In these maps  the main source of initial  
conditions sensitivity is given by the fact that nearby starting orbits can be separated by the  
discontinuities of the map. For this reason the initial condition sensitivity is estimated when we estimate how near we go   to the discontinuity points.  }.
 In the present work this argument is  extended.  
 The results we present  give an {\em upper}  bound to the local dimension of the measure at the point  $y$ in terms of the generalized recurrence indicator.  
Moreover, these generalized  recurrence indicators are also easy to be estimated numerically by looking to the behavior of a "typical" orbit, measuring its
 first entrance time in a decreasing sequence of balls centered in $y$. 
These results and the  ones given in \cite{BS} can be combined to have upper and lower bounds on the (upper and lower) local dimension of general systems. 
 One last remark is that the following results also hold  in systems with an infinite
invariant measure.

\section{Recurrence and dimension} 
 
In the following we will consider a discrete time dynamical system $(X,T)$ were $X$ is a separable metric space equipped with a Borel  locally finite measure  $\mu $  and $T:X\rightarrow X$ is a measurable map (we remark that we do not assume 
$\mu(X)=1$).

Let us consider the first entrance time of the orbit of $x$ in the ball $B(y,r)$ with center $y$ and radius $r$

$$\tau _r (x ,y)=\min ( \{n\in{\bf N},n>0,  T^n(x)\in B(y,r)\} ).$$ 

\noindent By this let us define the quantitative recurrence indicators 

$$\overline R(x,y) =\stackunder {r\rightarrow 0}{limsup}\frac {log(\tau _r (x,y))}{-log(r)}=
\stackunder {n\rightarrow \infty}{limsup}\frac {log(\tau _{2^{-n}} (x,y))}{n}$$ 

$$\underline R(x,y) =\stackunder {r\rightarrow 0}{liminf}\frac {log(\tau _r(x,y))}{-log(r)}= 
 \stackunder {n\rightarrow \infty}{limsup}\frac {log(\tau _{2^{-n}} (x,y))}{n}.$$
If for some $r$ $ \tau _r (x ,y)$  is infinite then $\overline R(x,y)$ and
$\underline R(x,y)$
are set to be equal to infinity.
 The indicators $\overline R(x)$ and $\underline R(x)$ of quantitative recurrence  defined in \cite{BS} are obtained as a special case, 
 $\overline R(x)=\overline R(x,x)$, $\underline R(x)=\underline R(x,x)$.

We state some first properties of $R(x,y)$. The proof follows directly from the definitions.

\begin{proposition}\label{inizz}$R(x,y)$ satisfies the following properties 
\begin{itemize}  
\item $\overline R(x,y) =\overline R(T(x),y)$, $\underline R(x,y) =\underline R(T(x),y)$.
 
\item If $T$ is Lipschitz, then  $\overline R(x,y)\geq\overline R(x,T(y))$, $\underline R(x,y)\geq\underline R(x,T(y))$ . 

\item If $T $ is $\alpha-Hoelder$, then  $\overline R(x,y)\geq \alpha \overline R(x,T(y))$, $\underline R(x,y)\geq\alpha \underline R(x,T(y))$.   
\end{itemize} 
\end{proposition}

Now we are interested to prove relations with dimension.
If $X$ is a metric space and $\mu$ is a measure on $X$ 
the upper local dimension at $x\in X $ is defined  as 

$$\overline {d}_\mu (x)  =\stackunder{r\rightarrow 0}{limsup}\frac {log(\mu (B(x,r))}{log(r)}=\stackunder {k\in {\bf N}, k\rightarrow \infty}{limsup} \frac{-log (\mu (B(x,2^{-k}))}k$$
 
 \noindent the lower local dimension $\underline {d}_\mu (x)$ is defined in an analogous way by replacing $limsup$ with $liminf$.

In general (even in examples that are interesting in dynamical system 
theory) $\underline {d}_\mu (x) $ and $\overline {d}_\mu (x) $ can  differ
on a positive measure set. 
If $\overline {d}_\mu (x)=\underline {d}_\mu (x)=d$ almost everywhere
the system is called exact dimensional. In this case all notions of dimension
of a measure (Hausdorff, box counting, information dimension) will coincide  (See for example the book \cite{P}) and then we have a precise description of the fractal structure of the system.
. For these and other reasons is important to have estimations for both $\underline {d}_\mu (x) $ and $\overline {d}_\mu (x) $.

With the above notations, Theorem 1  of \cite{BS} can be rewritten as follows

\begin{theorem}
If $X$ is a closed subset of $\bf{R}^n$ then 
for almost each $x\in X$ 
$$\overline{R}(x,x)\leq \overline{d}_\mu (x)\ , \ \underline{R}(x,x)\leq \underline{d}_\mu (x).$$
\end{theorem} 
In uniformly hyperbolic systems \cite{BS} also proved that recurrence and dimension are a.e. equal.
The equality also holds in some nonuniformly hyperbolic example, however
it is not difficult to see (\cite{BS} example 3) that there are uniquely ergodic irrational 
 rotations $(S^1,x\rightarrow x+\alpha)$ such that $\underline{R}(x,x)< \underline{d}_\mu (x)$ for each $x\in S^1$. 
In such systems and in general systems ${R}(x,x)$  then gives only a lower
bound for the dimension. We will see how it is possible to obtain a general upper bound for the dimension in term of $R(x,y)$.
 In \cite{BGI}  is indeed proved 

\begin{lemma}\label{Ganzo!} 
Let $(X,T)$ be as above, $\mu$ is an invariant measure, and $y\in X$. If $\alpha>{\underline {d}_\mu(y)}^{-1}$
then for $\mu$-almost all $x \in X$ it holds
$$\liminf_{n\to \infty}n^\alpha \,\min_{i}\,d(y,T^nx)=\infty.$$
\end{lemma}

Here we reformulate and extend this fact in the following way

\begin{theorem}\label{GAN} If $\mu $ is invariant,
for each fixed $y$ 
\begin{equation}\label{thm4} \underline {R}(x,y)\geq \underline {d}_{\mu}(y)\ , \ \overline {R}(x,y)\geq \overline {d}_{\mu}(y)\end{equation}
\noindent holds for $\mu$ almost each $x$.
\end{theorem}

{\em Proof.} 
First we prove $\underline {R}(x,y)\geq \underline {d}_{\mu}(y)$.
We remark that if $n^{-\alpha}\leq r \leq (n+1)^{-\alpha}$, since $\tau_r(x,y)$ is decreasing in $r$  then $$\frac {log(\tau_{n^{-\alpha}}(x,y))}{-log((n+1)^{-\alpha})}\leq \frac {log(\tau_r(x,y))}{-log(r)}\leq \frac {log(\tau_{(n+1)^{-\alpha}}(x,y))}{-log(n^{-\alpha})} $$ by this we can see that 
$\stackunder{r\rightarrow  0}{liminf} \frac {log(\tau_r(x,y))}{-log (r)}=
\stackunder{n \in {\bf N},n\rightarrow \infty}{liminf} \frac {log(\tau_{n^{-\alpha}}(x,y))}{-log (n^{-\alpha})}$.
Now Lemma \ref{Ganzo!} implies that if $n$ is big enough $\tau_{n^{-\alpha}}(x,y)\geq n $ for each $\alpha> \frac 1{\underline{d}_\mu(y)}$. Then  $\stackunder{n \in {\bf N},n\rightarrow \infty}{liminf} \frac {log(\tau_{n^{-\alpha}}(x,y))}{-log (n^{-\alpha})}\geq \frac 1\alpha .$
Since $\alpha$ can be chosen as near as we want
to $\frac 1{\underline{d}_\mu(y)}$ we have the statement.

Now we prove $\overline {R}(x,y)\geq \overline {d}_{\mu}(y)$. Suppose $d'< \overline{d}_\mu(y)$,
let us consider the set   $$A(d',y)=\{x\in X| \overline{R}(x,y)<d'\}.$$

\noindent By the assumption on the dimension, if $0<d'<d<\overline {d}_\mu (y) $ then there is a sequence $n_k$ such that
\begin{equation}\label{misu} \mu (B(y,2^{-n_k}))<2^{-dn_k}\ \mbox{ for each } k.
\end{equation}

On the other side  for each $x\in A(d',y)$ the relation $ \tau_{2^{-n}}(x,y)<2^{d'n}$ must hold eventually.
Let us consider $C(m)=\{x\in A(d',y)| \forall n\geq m, \tau_{2^{-n}}(x,y)<2^{d'n}\}$. This is an increasing sequence of sets. If we prove that $\stackunder {m \rightarrow \infty}{liminf}\mu(C(m))=0$
the statement is proved.
By the definition of $C(m) $ we see that $$C(n_k)\subset \stackunder {i\leq 2^{d'n_k}}{\cup} T^{-i}(B(y,2^{-n_k}))$$
\noindent the latter is made of $2^{d'n_k}$ sets, whose measure can be 
estimated by Eq. \ref{misu}, because $T$ is measure preserving.
Then $\mu (C(n_k))\leq 2^{d'n_k}*2^{-dn_k} $ and $\mu (C(n_k))$ goes to $0$
as $k\rightarrow \infty $.
${\Box} $

\begin{remark} If the measure $\mu $ is not invariant  inequality \ref{thm4} can fail
at some point. For example let us consider a system $(X,T)$, where
the map $T$ sends all the space $X$ in a   point $y$  ($\forall x, T(x)=y$) with  $\overline{d}_\mu (y)>0$.
Here $\overline{R}(x,y)< \overline{d}_\mu (y) $. We remark that in this example
the inequality fails only at one point ($y$). Next results shows that even
when the measure is not preserved the inequality can fail only on a zero measure set.
\end{remark}

In the previous result the invariance of the measure was an important ingredient. The following results (where $x$ is fixed and $y$ varies) are more general, they do not require the invariance of $\mu$. 
\begin{theorem}\label{2}
For each $x\in X$
$$\overline{R}(x,y)\geq \overline{d}_\mu (y)  \ , \ \underline{R}(x,y)\geq \underline{d}_\mu (y)$$\noindent for  $\mu$ almost each $y$.
\end{theorem}

\begin{theorem}\label{3} For each $x\in X$
the set $Y_h\subset X $ such that $$Y_h=\{y\in X , {\underline R}(x,y) \leq h\} $$ \noindent has
Hausdorff dimension $\leq h$. 
\end{theorem}

We remark that since obviously $\overline{R}(x,y)\geq \underline{R}(x,y)$
then the above result holds also with $\overline{R}(x,y)$ instead of $\underline{R}(x,y)$.

\begin{remark}\cite{KS}
Before to prove these results we remark that one cannot expect in general 
stronger results like $\overline{R}(x,y)\stackunder{\mu \ a.e.}= \overline{d}_\mu (y)$.
This can be realized by thinking to the following trivial example:
let us consider a periodic rotation  $(S^1, x\rightarrow x+\alpha, \lambda )$ with $\alpha \in {\bf Q}$ and $\lambda$  is the Lesbegue measure, here $\overline R(x,y)=\underline R(x,y)=\infty $ for each $y$ that in not contained in the orbit of $x$ (that is a finite set)
while $\lambda $ has dimension $1$. A less trivial example is in a certain sense a small perturbation of this latter one.

An irrational $\alpha$ is said to be of type $\nu$ if $$\nu=sup\{ \beta|\stackunder {n\rightarrow \infty} {liminf} j^{\beta}\stackunder {n\in {\bf N}}{min}|j\alpha -n|=0\}.$$ Lesbegue almost each irrational is of type $1$, but there are 
irrationals with type $>1$.
From the main result of \cite{KS} it can be deduced  (with some techical work) 
that an irrational rotation with angle of type $\nu>1$ satisfies $\overline{R}(x,y)=\nu>1$ almost everywhere (while $\underline{R}(x,y)=1$ a.e.).

\end{remark}

 \section{ Proof of Theorems  \ref{2} and \ref{3}}

Theorems \ref{2} and \ref{3} come  from the following more general results. Let us  consider a sequence $x_i:{\bf N}\rightarrow X$, we define  recurrence indicators indicating how the sequence  
 comes near some given points. 
 For this let us consider $y\in X$, and the first entrance time of $x_i$ in a ball with center $y$ 

$$\tau (x_i ,y,r)=min\{n\in{\bf N}, x_n\in B(y,r)\}.$$ 

\noindent Let us define the quantitative recurrence indicators 

$$\overline R(x_i,y) =\stackunder {r\rightarrow 0}{limsup}\frac {log(\tau(x_i,y,r))}{-log(r)} =\stackunder {n\rightarrow \infty }{limsup}\frac {log(\tau(x_i,y,2^{-n}))}{n}$$ 

$$\underline R(x_i,y) =\stackunder {r\rightarrow 0}{liminf}\frac {log(\tau(x_i,y,r))}{-log(r)}=\stackunder {n\rightarrow \infty }{liminf}\frac {log(\tau(x_i,y,2^{-n}))}{n} $$

Theorem \ref{3} comes from the following proposition
\begin{proposition} For each sequence $x_i$
the set $Y_h\subset X $ such that $Y_h=\{y\in X , {\underline R}(x_i,y) \leq h\} $ has
Hausdorff dimension $\leq h$. 
\end{proposition}

{\em Proof.}
We have that $\forall y \in Y_h$ $ \stackunder {k \rightarrow \infty} {minlim} \frac{log(\tau (x_i,y,2^{-k}))}k \leq h$.
This means that $\forall \epsilon>0, \forall y \in Y_h$ and $\forall k_0 \in {\bf N}$ there is  
$k>k_0$ and an index   $j$ with $j \leq 2^{(h+\epsilon) k}$ with
$y\in B(x_j,2^{-k})$.
 
Let us call $S_{\epsilon, k}$ the union of all the balls $B(x_j,2^{-k})$ for all index
$j$ such that
 $j \leq 2^{(h+\epsilon)k}$; $S_{\epsilon, k}=\stackunder {j \leq 2^{(h+\epsilon)k}}{\cup}B(x_j,2^{-k})$.
 Then for each $k_0$
$$ Y_h\subset \stackunder {k\geq k_0}\cup S_{\epsilon, k}     $$

\noindent by this  $Y_h $ (and each $ S_{\epsilon, k},k>k_0$) is covered by a family of
balls of diameter less than $2^{-k_0}$ and  we can  estimate the $d-$dimensional Hausdorff measure of $ S_{\epsilon, k}$ 

$$ {\cal H}_{2^{-k+1}}^d ( S_{\epsilon, k})\leq 2^{(h+\epsilon)k+1}2^{(-k+1)d}
= 2^{k(h+\epsilon -d)+1+d}$$
and  $$ {\cal H}_{2^{-k_0}} ^d (Y_h)\leq \stackunder{k\leq k_0}{\sum} 2^{1+d}2^{k(h+\epsilon -d)}$$
if $d>h+\epsilon$ we can set $k_0$ so big that $  {\cal H}_{2^{-k_0}} ^d(Y_h)
\leq \delta $ for each fixed $\delta$ and $Y_h$ is covered by balls of arbitrary small
size.
This proves $ {\cal H}^d(Y_h)=\stackunder {k_0\rightarrow
\infty}{lim}{\cal H}_{2^{-k_0}} ^d(X_h)=0 $ for each   $ d>h+\epsilon$.
Since $\epsilon $ is arbitrary the statement follows.$\Box $

\begin{remark}By \cite{BS} (example 3) we have that if $\alpha $ is ``well approximable'' by rationals then $\underline R(x,x)<1$. By theorem \ref{3} the set of other points $y$ such  that  $\underline R(x,y)<d<1$ is very small, indeed it must have dimension less or equal than $d$.
\end{remark}

Theorem \ref{2} comes from
\begin{proposition}\label{propo}
For $\mu $ almost each $y\in X$
$$\overline{R}(x_i,y)\geq \overline{d}_{\mu}(y) \ ,\  \underline{R}(x_i,y)\geq \underline {d}_{\mu}(y).$$ 
\end{proposition}

The proof of proposition \ref{propo}  is based on the following lemmas

\begin{lemma} Let  $A=\{y\in X,\overline{d}_\mu(y)> d \}$.  If $h<d$  and $$Y_h=\{y\in A, s.t.\ \overline{{ R}}(x_i,y)<h\}$$                          
\noindent then  $\mu (Y_h)=0$.  
 \end{lemma} 

{\em Proof. }Let us  consider $0<\varepsilon < h-d$ and 
\[
Y^{n}_{h}=\{y\in A\, s.t.\, \forall m>n\, log(\tau (x_i,y,2^{-m}))<(h+\varepsilon)m \}\]
 we have that \( Y_{h}\subset \stackunder {n\geq n_0}{\cup} Y_{h}^{n} \) , \( Y^{n}_{h}\subset Y^{n+1}_{h}. \)
If we prove that \( \mu (Y_{h}^{n})=0 \) eventually with respect to $n$ the assertion is proved.

If \( y\in Y^{n}_{h} \) then \( \forall m>n\,\exists x, s.t. \  y\in B(x_i,2^{-m}) \) where \( i<2^{m(h+\varepsilon )} \) in other
words if we consider the set of all ball of radius \( 2^{-m} \) with centers $x_i$ with \( i<2^{m(h+\varepsilon )} \)

\[
{\cal B} ^{m}=\{B(x_i,2^{-m})\, s.t.\, i<2^{(m(h+\varepsilon ))}\}\]
 we have that \( \forall m>n\  , Y_{h}^{n}\subset \stackunder{\beta \in {\cal B} ^{m}}{\bigcup }\beta  \).

For each \( y \) \( \in A \) we have that 
$\overline{d}(y)= \stackunder {n\rightarrow \infty } {\limsup} \frac{-\log \mu (B(y,2^{-n}))}{n}> d$,
 this implies that \( \forall y\in A \) there exist an infinite sequence \( B(y,2^{-n_{k}}) \)
of balls centered in \( y \) ,with radius \( 2^{-n_{k}d} \) such that \( \mu (B(y,2^{-n_{k}}))\leq 2^{-n_{k}d} \)
. Let us call this family of balls \( y-\)estimated balls.

Now let us consider the balls in \( {\cal B} ^{n} \) for which we have an estimation about their measure: we say that a ball in \( \beta \in {\cal B} ^{n} \) is ``nice'' if there exist an \( y\textrm{ } \)such that \( \beta  \) is contained in
some \( y\)-estimated ball of radius \( 2^{-n+1} \)
found above (we recall that all the
balls in \( {\cal B} ^{n} \) have radius \( 2^{-n} \)), thus if \( \beta  \)
is nice then \( \mu (\beta )\leq 2^{-(n-1)d} \). Every \( y\in Y_{h}^{n} \)
has a sequence of \( y-\)estimated balls, let us consider one of these balls
\( \beta (y,2^{-k+1}) \): \( y \) is also contained in a ball \( \beta '\in {\cal B} ^{k} \),then
\( \beta '\subset \beta(y,2^{-k+1}) \). This implies that \( \forall j\geq n \)
each point of  \( Y_{h}^{n} \) is contained in some ``nice'' ball with radius
not greater than \( 2^{-j} \) that is: \(  Y_{h}^{n}\subseteq \stackunder{m\geq n}{\bigcup}\stackunder {\beta \in \{nice\, balls\, in\, {\cal B} ^{m}\}}{\bigcup}\beta  \).
Now we are ready to estimate the total measure of the nice balls: we remark
that the number of nice balls with radius \( 2^{-m}\textrm{ } \)is not greater
than \( 2^{m(h+\varepsilon )+1} \) and the measure of a nice ball is not greater
than the measure of the corresponding \( y- \) estimated ball. This implies
that \( \forall j>n \)$$ \mu (Y_{h}^{n})\leq \stackunder {m\geq j}{\sum}2^{m(h+\varepsilon )+1}2^{-(m-1)d} $$ if $n$ is big enough this sum can be set as small as we want, then \( \mu (X^{n}_{h})=0 \) . {$\Box $}

\begin{lemma}\label{aboveo}
Let $d,c,\delta >0$, let $B=\{y\in X,  \underline {d}_{\mu} (y)>d+c \}$. If $A=\{y\in B   \ s.t. \ \underline{R}(x_i,y) \leq d-\delta \}$, then $\mu(A)=0$. 
\end{lemma}

{\em Proof. } Conversely let us suppose that $\mu(A)>0$.
Since $\forall x\in  A $, $\underline {d}_{\mu}(x) \geq d+c $ then $ x\in A $, implies that  if $m $ is big enough (depending on $x$)  \begin{equation} 
\label{supra}\mu(B(x,2^{-m}))<2^{-m(d+\frac c2 )}\end{equation}
\noindent then there is an $\overline{m}>0$ and a set $A'\subset A$ with $\mu (A')>0$ such that if $x\in A'$ $\forall m >\overline {m} \ \mu(B(x,2^{-m}))<2^{-m (d+\frac c 2)}$ {\em uniformly} on all $A'$.

By the definition of $A'$ for each $k$, each $x\in A'$ is contained in some ball
$B(x_j,2^{-i})$ with $j>k$ and $j<2^{di}$, that is 
 $\forall k\in {\bf N}$ $A\subset \stackunder { \{ j|k \leq j\} } {\cup} B(x_j,2^{-\frac {log(j)}{d}})$.

Now the measure of these balls can be estimated as before  by Eq. \ref{supra} and then the total measure of $A'$ can be estimated as in the previous proof, concluding
that $\mu(A)=0$.$\Box$

{\em  Proof of proposition \ref{propo} } If   conversely  $\overline{R}(x,y)< \overline{d}(x)$
on a set $A'$ with $\mu(A')>0$ we can find a constant $c$ and a set $A''$, $\mu(A'')>0$ such that $\overline{R}(x,y)<c< \overline{d}(x)$ on $A''$, by lemma  16 we obtain  $\mu(A'')=0$.
Similarly the other inequality can be obtained $\Box $

\end{document}